\documentclass{article}

\usepackage{amsfonts}
\usepackage{amssymb}
\usepackage{amsmath}
\usepackage{amsthm}

\theoremstyle{plain}
\newtheorem{theorem}{Theorem}[section]

\theoremstyle{definition}

\newtheorem{example}[theorem]{Example}

\newcommand{\NN}{{\mathbb N}}

\usepackage[T1]{fontenc}
\usepackage[adobe-utopia]{mathdesign}
\usepackage{color}

\title{Transfunctions}
\author{Piotr Mikusi\'nski\\
Department of Mathematics\\
University of Central Florida\\
 Orlando, Florida, USA}

\begin{document}
\maketitle

\begin{abstract}  Maps between spaces of measures on measurable spaces $(X,\Sigma_X)$ and $(Y, \Sigma_Y)$ are treated as generalized functions between $X$ and $Y$.  
\end{abstract}

\bigskip



\section{Introduction}

The defining property of a function is that it has values at points.  More precisely, if $f$ is a functions from $X$ into $Y$, then for every $x\in X$, $f(x)\in Y$. While this is a very important and useful idea in mathematics, there are many situations in mathematics and in applications of mathematics where the definition of a function is too restrictive. There are many different generalizations of the idea of a function, including multivalued functions, functions defined almost everywhere, fuzzy functions, Mikusi\'nski operators, Schwartz distributions, hyperfunctions, ultradistributions, Colombeau generalized functions, or Bohemians. They are motivated by different applications and are constructed with a variety of mathematical tools.

In this short note we propose a new framework for generalizing functions. The motivation for our definition is the fact that in real world the input is rarely a precise element of the domain and the same is true about the output.  We prefer to think of the input as a ``set of points with probabilities'' and the output as a ``set of outputs with probabilities''. More precisely, we are mapping measures in the domain space to measures in the target space.  For mathematical convenience, we are not restricting the measures in the domain and range to probability measures, but instead we consider arbitrary finite measures. We call these objects {\it transfunctions}.

While the intuition behind transfunctions is similar to that of fuzzy functions, the    
mathematical formalisms of these two approaches are very different.  

This short note should be treated as an introduction of transfunctions.  We only present here the general definition and some examples.

\section{The general definition}

Let $(X,\Sigma_X)$ and $(Y, \Sigma_Y)$ be measurable spaces, with $\Sigma_X$ and  $\Sigma_Y$ the respective $\sigma$-algebras, and $\mathcal{M}_X$ the set of all finite measures on $\Sigma_X$, and $\mathcal{M}_Y$ the set of all finite measures on $\Sigma_Y$.   
By a transfunction from $(X,\Sigma_X,\mathcal{M})$ to $(Y, \Sigma_Y,\mathcal{N})$ we mean a map
$$
\Phi : \mathcal{M}  \to \mathcal{N},
$$
where $\mathcal{M}\subset \mathcal{M}_X$ and $\mathcal{N}\subset\mathcal{M}_Y$ are positive cones.
We consider the following properties of transfunctions:
\begin{description}
\item[Weakly additive:] $\Phi(\mu_1+\mu_2)= \Phi(\mu_1)+\Phi(\mu_2)$ for all $\mu_1,\mu_2\in\mathcal{M}$ such that $\mu_1 \perp \mu_2$,
\item[Strongly additive:] $\Phi(\mu_1+\mu_2)= \Phi(\mu_1)+\Phi(\mu_2)$ for all $\mu_1,\mu_2\in\mathcal{M}$,
\item[Homogeneous:] $\Phi(\alpha \mu)=\alpha \Phi(\mu)$ for any $\alpha > 0$,
\item[Monotone:] $\Phi(\mu_1) \leq \Phi(\mu_2)$, if $\mu_1\leq \mu_2$,
\item[Measure preserving:] $\|\Phi(\mu)\|=\|\mu\|$,
\item[Bounded:] $\|\Phi(\mu)\|\leq C\|\mu\|$ for some $C>0$ and all $\mu\in \mathcal{M}$,
\item[Continuous:] $\Phi(\mu_n) \to \Phi(\mu)$ whenever $\mu_n \to \mu$.
\end{description}
Since different types of convergence of measures can be considered, they have to be specified when discussing continuity of transfunctions.

\section{Examples}

\begin{example}\label{function}
If $(X,\Sigma_X)$ and $(Y, \Sigma_Y)$ are measurable spaces and $f:(X,\Sigma_X)\to (Y, \Sigma_Y)$ is a measurable function, then $\Phi : \mathcal{M}_X  \to \mathcal{M}_Y$ defined by
$$
\Phi_f(\mu)(S)=\mu(f^{-1}(S)) 
$$
is a strongly additive, homogeneous, monotone, and measure preserving transfunction. Moreover, $\Phi_f$ is continuous with respect to strong convergence in $\mathcal{M}_X$ and $\mathcal{M}_Y$.
\end{example}

\begin{example}
Let $X=\{x_1,\dots , x_m\}$, $Y=\{y_1,\dots ,y_n\}$, $\Sigma_X=2^X$, $\Sigma_Y=2^Y$. Then
$\mathcal{M}_X$ can be identified with $[0,\infty)^m$ and $\mathcal{M}_Y$ can be identified with  $[0,\infty)^n$, and a transfunction $\Phi : \mathcal{M}_X  \to \mathcal{M}_Y$ can be identified with a map
$$
\Phi : [0,\infty)^m \to [0,\infty)^n.
$$

If $A=(a_{i,j})$ is an $n\times m$ matrix with nonnegative entries, then the transfunction defined by
$$
\Phi (\mu)= A \mu
$$
is strongly additive, homogeneous, monotone, and bounded. If $\sum_{i=1}^n a_{i,j}=1$ for every $j\in\{1,\dots ,m\}$, then $\Phi$ is measure preserving. Note that $\Phi$ corresponds to a function $f:X\to Y$ (as defined in Example \ref{function}) if an only if $A$ is a 0-1 matrix and for every $i\in\{1,\dots ,m\}$ there is exactly one $j\in\{1,\dots ,n\}$ such that $a_{i,j}\neq 0$.
\end{example}

The above example can be generalized to infinite sets.

\begin{example}
Let $X=\{x_1,x_2,\dots \}$, $Y=\{y_1,y_2,\dots \}$, $\Sigma_X=2^X$, and $\Sigma_Y=2^Y$. Then
$$
\mathcal{M}_X = \mathcal{M}_Y = l_+^1= \{(\varepsilon_n)\in l^1: \varepsilon_n\geq 0\}.
$$

If $A=(a_{m,n})$ is an infinite matrix with nonnegative entries such that 
$$
\sum_{m=1}^\infty a_{m,n} < M
$$
 for some constant $M$ and all $n\in\NN$, then the transfunction defined by
$$
\Phi (\mu)= A \mu
$$
is strongly additive, homogeneous, and monotone. If $\sum_{m=1}^\infty a_{i,j}=1$ for every $n\in\NN$, then $\Phi$ is measure preserving. 
\end{example}

\begin{example}
 Let $(X,\Sigma_X)$ be a measurable space, let $\rho$ be a finite positive measure on $(Y,\Sigma_Y,\rho)$, and let $\varphi \in L^\infty (X\times Y, \Sigma_X \times \Sigma_Y)$ be nonnegative. For $\mu\in\mathcal{M}_X$ we define
$$
\Phi(\mu)(B)=\int_{X\times B} \varphi \, d(\mu \times \rho),
$$
where $B\in \Sigma_Y$. Clearly, $\Phi$ is strongly additive, homogeneous, and monotone. Since,
$$
\|\Phi(\mu)\|=\int_{X\times Y} \varphi \, d(\mu \times \rho)\leq \|\varphi\|_\infty \|\rho\| \|\mu\|,
$$
$\Phi$ is bounded.
\end{example}

\begin{example} Let $\Phi$ be a transfunction from $(X,\Sigma_X,\mathcal{M})$ to $(Y,\Sigma_Y,\mathcal{N})$ and let $f:Y\to [0,\infty)$ be a multiplier on $\mathcal{N}$, that is $f\nu\in \mathcal{N}$ whenever $\nu \in \mathcal{N}$. Then 
$$
(f\Phi)(\mu)(B) = f(\Phi(\mu))
$$
is transfunctions from $(X,\Sigma_X,\mathcal{M})$ to $(Y,\Sigma_Y,\mathcal{N})$.
\end{example}

\begin{example} Let $\Phi$ be a transfunction from $(X,\Sigma_X,\mathcal{M})$ to $(Y,\Sigma_Y,\mathcal{N})$ and let $g:X\to [0,\infty)$ be a multiplier on $\mathcal{M}$. Then 
$$
(\Phi g)(\mu) = \Phi(g\mu)
$$
is transfunctions from $(X,\Sigma_X,\mathcal{M})$ to $(Y,\Sigma_Y,\mathcal{N})$.
\end{example}

\begin{example}
 Let $\Phi$ be a transfunction from $(X,\Sigma_X,\mathcal{M})$ to $(Y,\Sigma_Y,\mathcal{N})$ and let $\rho$ be a finite positive measure on $(Y,\Sigma_Y,\rho)$. Then 
$$
\Psi(\mu) = \max \{\Phi(\mu),\rho\}
$$
is transfunctions from $(X,\Sigma_X,\mathcal{M})$ to $(Y,\Sigma_Y,\mathcal{N})$.
\end{example}

\begin{example} If $\mu$ is a measure on some $\sigma$-algebra $\Sigma$ on a set $Z$ and $S\in\Sigma$, then by $\pi_S \mu$ we will denote the projection of $\mu$ onto $S$, that is, $\pi_S\mu(A)=\mu(A\cap S)$.

Let $\Phi$ be a transfunction from $(X,\Sigma_X,\mathcal{M})$ to $(Y,\Sigma_Y,\mathcal{N})$. If $A\in\Sigma_X$ and $B\in\Sigma_Y$, then 
$$
\Phi\pi_A \quad \text{and} \quad \pi_B\Phi
$$
are transfunctions from $(X,\Sigma_X,\mathcal{M})$ to $(Y,\Sigma_Y,\mathcal{N})$.
\end{example}

\begin{example}
 Assume that $(Y,\odot)$ is a semigroup and that $\odot : Y \times Y \to Y$ is a measurable function with repect to the product $\sigma$-algebra generated by $\Sigma_Y\times\Sigma_Y$.

Let $\Phi$ and $\Psi$ be transfunctions from $(X,\Sigma_X,\mathcal{M})$ to $(Y, \Sigma_Y,\mathcal{N})$. Then
$$
(\Phi \odot \Psi) (\mu)(B)= (\Phi(\mu)\times\Psi(\mu)) (\odot^{-1}(B)).
$$
defines a transfunction from $(X,\Sigma_X,\mathcal{M})$ to $(Y, \Sigma_Y,\mathcal{N})$.
\end{example}

\section{Remarks}

We are interested in studying the structure of transfunctions, their algebraic and topological properties, as well as potential applications to problems where point to point functions are inadequate. We would also like to better understand how transfunctions are related to other  generalizations of functions.

\end{document}